\documentclass[1 [leqno,11pt]{amsart}
\usepackage{amssymb, amsmath}
\def\refer#1{~\ref{#1}}
\def\refeq#1{~(\ref{#1})}
\def\ccite#1{~\cite{#1}}

\def\longformule#1#2{
\displaylines{ \qquad{#1} \hfill\cr \hfill {#2} \qquad\cr } }
\def\inte#1{
\displaystyle\mathop{#1\kern0pt}^\circ }

\def\supetage#1#2{
\sup_{\scriptstyle {#1}\atop\scriptstyle {#2}} }

\let\al=\alpha

\let\e=\varepsilon

\let\lam=\lambda

\def\buh{{v}^h}

\let\p=\psi

\let\D=\Delta

\let\wt=\widetilde

\def\with{\quad\hbox{with}\quad}
\def\andf{\quad\hbox{and}\quad}

\def\virgp{\raise 2pt\hbox{,}}
\def\cdotpv{\raise 2pt\hbox{;}}

\def\eqdefa{\buildrel\hbox{\footnotesize def}\over =}

\def\C{{\mathop{\bf C\kern 0pt}\nolimits}}
\def\DD{{\mathop{\bf D\kern 0pt}\nolimits}}
\def\K{{\mathop{\bf K\kern 0pt}\nolimits}}
\def\N{{\mathop{\mathbb N\kern 0pt}\nolimits}}
\def\PP{{\mathop{\mathbb P\kern 0pt}\nolimits}}
\def\Q{{\mathop{\mathbb Q\kern 0pt}\nolimits}}
\def\R{{\mathop{\mathbb R\kern 0pt}\nolimits}}
\def\SS{{\mathop{\mathbb S\kern 0pt}\nolimits}}
\def\ZZ{{\mathop{\mathbb Z\kern 0pt}\nolimits}}
\def\TT{{\mathop{\mathbb T\kern 0pt}\nolimits}}
\def\BB{{\mathop{\mathbb B\kern 0pt}\nolimits}}
\def\PP{{\mathop{\mathbb P\kern 0pt}\nolimits}}

\newcommand{\ds}{\displaystyle}

\def\dive{\mathop{\rm div}\nolimits}

\def\no{\noindent}
\def\na{\nabla}
\def\p{\partial}

\newcommand{\beq}{\begin{equation}}
\newcommand{\eeq}{\end{equation}}
\newcommand{\ben}{\begin{eqnarray}}
\newcommand{\een}{\end{eqnarray}}
\newcommand{\beno}{\begin{eqnarray*}}
\newcommand{\eeno}{\end{eqnarray*}}

\newtheorem{defi}{Definition}[section]
\newtheorem{thm}{Theorem}
\newtheorem{lem}[defi]{Lemma}
\newtheorem{rmk}[defi]{Remark}

\newtheorem{prop}[defi]{Proposition}

\begin{document}

\title[Sums of large global solutions] {Sums of large global
solutions \\to the incompressible  Navier-Stokes equations} {}

\author[J.-Y. CHEMIN]{Jean-Yves Chemin}
\address [J.-Y. Chemin]%
{Laboratoire J.-L. Lions, UMR 7598 \\
Universit\'e Pierre et Marie Curie, 75230 Paris Cedex 05, FRANCE }
\email{chemin@ann.jussieu.fr}
\author[I. GALLAGHER]{Isabelle GALLAGHER}
\address [I. Gallagher]%
{Institut de Math\'ematiques de Jussieu \\
Universit\'e Paris Diderot, 75251 Paris Cedex 05, FRANCE} \email
{gallagher@math.jussieu.fr}
\author[P. Zhang]{Ping Zhang}%
\address[P. Zhang]
{Academy of
Mathematics $\&$ Systems Science, CAS\\
Beijing 100080, CHINA\\
and  Hua Loo-Keng Key Laboratory of Mathematics, Chinese Academy of
Sciences.} \email{zp@amss.ac.cn}

\maketitle

\begin{abstract}
Let~${\mathcal G}$ be the (open) set of~$\dot H^{\frac 1 2}(\R^3)$
divergence free vector fields generating  global smooth solutions
 to the three
dimensional incompressible Navier-Stokes equations. We prove that
any element of~${\mathcal G}$ can be perturbed by an arbitrarily
large, smooth divergence free vector field which varies slowly in
one direction, and the resulting vector field (which remains
arbitrarily large) is an element of~${\mathcal G}$ if the variation
is slow enough. This result implies that through any point
in~${\mathcal G}$ passes an uncountable number of arbitrarily long
segments included in~${\mathcal G}$.\end{abstract}

\setcounter{equation}{0}
\section{Introduction}
\subsection{Setting of the problem and statement of the result}

Let us first recall the classical Navier-Stokes system for
incompressible fluids in three space dimensions:
$$ (NS)\qquad\left\{
\begin{array}{ll}
\p_tu +u \cdot\na u -\Delta u =-\na p \\
\dive u =0\\
u _{|t=0}=u_0
\end{array}\right.
$$
where $u (t,x)$ denotes the fluid velocity and $p (t,x)$ the
pressure. In this paper the space variable~$x$ is chosen in~$\R^3$.

All the solutions we are going to consider here are at least continuous in time with values
 in the Sobolev space~$\dot H^{\frac 1 2}(\R^3)$. It is well known that in that case, all
  concepts of solutions coincide and in particular we shall deal with
  $''$mild$"$
  solutions  of $(NS)$ (see for instance\ccite{lemarie1}).

In order to  specify the concept of large data, let us recall the
history of results concerning small data.  The first one states that
if the initial data~$u_0$ is such that~$\|u_0\|_{L^2} \|\nabla
u_0\|_{L^2} $ is small enough, (NS) has a global regular solution;
this was proved by J. Leray in his seminal paper\ccite{lerayns}.
Then, starting with the paper by H. Fujita and T. Kato
(see\ccite{fujitakato}), the following approach was developped: let
us denote by~$\BB$ the bilinear operator defined by
$$
\left\{\begin{array}{c} \ds\partial_t \BB (v,w) -\D \BB (v,w) =
\frac 12 \PP \dive (v\otimes w+
w\otimes v)\\
\BB (v,w) _{|t=0} =0
\end{array}
\right.
$$
where~$\PP$ denotes the Leray projection onto divergence free vector
fields. Then, it is easily checked that~$u$ is a solution of~$(NS)$
if and only if
$$
u=e^{t\D} u_0+\BB(u,u)
$$
which is something like Duhamel's formula. Then the theory of small
initial data reduces to finding a Banach space~$X$ of time-dependent
divergence free vector fields on~$\R^+\times\R^3$ such that~$\BB$ is
a bilinear map from~$X\times X$ to~$X$. An elementary abstract fixed
point theorem claims that if~$X$ is a
Banach space of time-dependent divergence free vector fields
on~$\R^+\times\R^3$ such that
$$
\|\BB(v,w)\|_X\leq C\|v\|_X \|w\|_X
$$
($X$ will be called from now on an adapted space), a solution
of~$(NS)$ exists in~$X$ and is global as soon as
$$
\|e^{t\D} u_0\|_X \leq (4C)^{-1}.
$$
The search of the largest possible adapted space~$X$ is a long story. It
started in~1964 with the paper\ccite{fujitakato} where the space~$X$
is defined by the norm
$$
\|u\|_X\eqdefa \sup_{t\geq 0} t^{\frac 1 4}\|\nabla u(t)\|_{L^2}.
$$
This corresponds to an initial data small in~$\dot H^\frac12(\R^3)$,
and it is shown in particular
 that the solution belongs to~$C(\R^+,\dot H^\frac12(\R^3)) \cap L^2(\R^+,\dot H^\frac32(\R^3))$.
After a number of important steps (see in
particular\ccite{kato},\ccite{gigamiyakawa},\ccite{weissler} and
\ccite{cannonemeyerplanchon}), the problem of finding the largest
adapted space was achieved by H. Koch and D. Tataru. They proved
in\ccite{kochtataru} that the space of time-dependent divergence
free vector fields on~$\R^+\times\R^3$ such that
$$
\|u\|_{X_{KT}} \eqdefa \sup_{t\geq 0} t^{\frac 1 2}
\|u(t)\|_{L^\infty} +\supetage{x\in \R^3} {R>0} \frac 1 {R^{\frac 3
2}} \left(\int_{P(x,R)} | u(t,y) | ^2dydt\right)^\frac12 <\infty
$$
where~$P(x,R)$ is the parabolic ball~$[0,R^2]\times B(x,R)$, is an
adapted space.

Now let us observe that the incompressible Navier-Stokes system is
translation and scaling invariant: if~$u$ is
a solution of~$(NS)$ on~$[0,T]\times\R^3$ then, for any
positive~$\lam$ and for any $x_0$ in~$\R^3$, the vector
field~$u_{\lam,x_0}$ defined by
$$
u_{\lam,x_0}(t,x) \eqdefa \lam u(\lam^2t,\lam (x-x_0))
$$
is also a solution of~$(NS)$
on~$[0,\lam^{-2}T]\times\R^3$. Thus, an adapted space must be
translation and scaling invariant in the following sense: a
constant~$C$ exists such that, for any positive~$\lam$ and for any
$x_0$ in~$\R^3$,
$$
C^{-1} \|u\|_X \leq \|u_{\lam,x_0}\|_X \leq C\|u\|_X.
$$
The second term appearing in the norm~$\|\cdot\|_{X_{KT}}$ above
comes from the fact that the solution of~$(NS)$ should be locally
in~$L^2$ in order to be able to define the product as a
locally~$L^1$ function. The relevant norm on the initial data
is~$\|e^{t\D} u_0\|_X$. In the case of the Koch and Tataru theorem,
this norm turns out to be equivalent to the norm of $\partial BMO$,
the space~
 of derivatives of $BMO$ functions. Of course, the
space of initial data which measures the size of the initial data
must be translation and scaling invariant. A remark due to Y. Meyer
(see~\cite{meyer}) is that the norm in such a space is always
greater than the norm in the Besov space~$\dot
B^{-1}_{\infty,\infty}$ defined by
$$
\|u\|_{\dot B^{-1}_{\infty,\infty}} \eqdefa \sup_{t\geq 0} t^{\frac
1 2} \|e^{t\D}u\|_{L^\infty}.
$$

This leads to the following definition of a large initial data for
the incompressible Navier-Stokes equations.
\begin{defi}
\label{definlargedata} {\sl  A divergence free vector field~$u_0$ is a
large initial data for the incompressible Navier-Stokes system if
its~$\dot B^{-1} _{\infty,\infty}$ norm is large.
}
\end{defi}

Let us point out that this approach using Duhamel's formula does not
use the very special structure of the incompressible Navier-Stokes
system. A family of results does use the special structure
of~$(NS)$: in those cases some geometrical invariance on the initial
data is preserved by the flow of $(NS)$ and this leads to some
unexpected conservation of quantities, which makes the problem subcritical and thus prevents  blow up. We refer for instance to \cite{ladyzhenskaya},
\cite{mahalovtitileibovich}, \cite{poncerackesideristiti}, or
\cite{ukhoiudo}, where special symmetries (like helicoidal, or
axisymmetric without swirl) allow to prove global wellposedness for
any data.

Some years ago, the first two authors investigated the possible
existence of large initial data (in the sense of
Definition\refer{definlargedata}) which have no preserved
geometrical invariance and which nevertheless generate global
regular solutions to $(NS).$ The first result in this direction was
proved in\ccite{cg} where such a family of large initial data was
constructed, with strong oscillations in one direction. The main
point of the proof is that for any element of this family, the first
iterate~$\BB(e^{t\D}u_0,e^{t\D} u_0)$ is exponentially small with
respect to the large initial data~$u_0$ in some appropriate norm.
Let us notice that this result does use the fine structure of the
non linear term of~$(NS)$: M. Paicu and the second author proved
in\ccite{gp} that for a modified incompressible Navier-Stokes
system, this family of initial data generates solutions that blow up
at finite time.

In\ccite{cg3}, the first two authors constructed another class of examples, in which the initial data has slow
variations in one direction. The proof of global regularity uses the fact that the 2D
Navier-Stokes equations are globally wellposed. The initial data
presented in the next theorem will be referred to in the following
as ``quasi-2D").
\begin{thm}[\cite{cg3}]
\label{theoquasi2D} {\sl  Let~$v_0^h=(v_0^1,v_0^2)$ be a two component,
smooth divergence free vector field on~$\R^3$ (i.e.~$v_0^h$ is
in~$L^2(\R^3)$ as well as all its derivatives), belonging, as well
as all its derivatives, to~$L^2(\R_{x_3}; \dot H^{-1}(\R^2))$; let
$w_0 = (w_0^h,w_0^3)$ be a three component, smooth divergence free vector field on~$\R^3$. Then there
exists a positive~$\e_0$ such that if~$\e\leq \e_0$, the initial
data
$$
u_{0,\e}(x) \eqdefa ( v_0^h+\e w^h_0,w_0^3 )(x_1,x_2 ,\e x_3)
$$
generates a unique, global solution~$u^{\e}$ of~$(NS)$. }
\end{thm}
\begin{rmk}
{\rm It is clear from the proof of\ccite{cg3} that the dependence of the
parameter~$\e_0$ on the profiles~$v_0^h$ and~$w_0$ is only through
their norms.}
\end{rmk}
Note that such an initial data may be arbitrarily large in the sense
of Definition\refer{definlargedata} (see~\cite{cg3}). We recall for the convenience of the reader the result
proved in\ccite{cg3}.
\begin{prop}[\cite{cg3}]
\label{sizedata}
{\sl  Let $(f,g)$ be in~$ {\mathcal S}(\R^2)\times {\mathcal S}(\R)$ and
define~$h_\e (x_h,x_3) \eqdefa f(x_h)g(\e x_3)$. We have, if $\e$ is
small enough,
$$
\|h_\e\|_{\dot B^{-1}_{\infty,\infty}(\R^3)} \geq \frac 1 4
\|f\|_{\dot B^{-1}_{\infty,\infty}(\R^2)} \|g\|_{L^\infty(\R)}.
$$
}
\end{prop}

\medskip

In this paper we consider the global wellposedness of the
Navier-Stokes equations with data which is the sum of an initial
data (which may be large) giving rise to a global solution, and a
quasi-2D initial data as presented above (which may also be large).
The theorem is the following.
\begin{thm}
\label{th1.1}
{\sl Let $u_0$, $v^h_0$ and~$w_0$ be three smooth
divergence free vector fields defined on~$\R^3$, satisfying

\begin{itemize}

\item $u_0$ belongs to~$\dot H^{\frac 1 2} (\R^3)$ and generates a
unique global solution to the Navier-Stokes equations;

\item $v_0^h=(v_0^1,v_0^2)$ is a horizontal  vector field on~$\R^3$ belonging, as well as all its
derivatives, to the space~$L^2(\R_{x_3}; \dot H^{-1}(\R^2))$;
\item $v^h_0(x_1,x_2,0)= w_0^3(x_1,x_2,0) =0$ for all $(x_1,x_2)\in\R^2.$

\end{itemize}
Then there exists a positive number~$\e_0 $ depending on~$u_0$ and
on norms of~$v_0^h$ and~$w_0$ such that for any $\e\in (0,\e_0]$,
there is a unique, global solution to the Navier-Stokes equations
with initial data
$$
u_{0,\e}(x) \eqdefa u_0(x)+(v^h_0 + \e w_0^h ,w_0^3)(x_1,x_2,\e
x_3).
$$
}
\end{thm}
\begin{rmk}
\label{linkwithabstract} {\rm Let~$u_0$ be any element of the (open)
set~${\mathcal G}$ of~$\dot H^{\frac 1 2}(\R^3)$ divergence free
vector fields generating global smooth solutions to~$(NS)$, and
let~$N$ be an arbitrarily large number. Then for any smooth
divergence free vector field~$f^h$ (over~$\R^2$) and scalar
function~$g$ (over~$\R$) satisfying~$\|f^h\|_{\dot
B^{-1}_{\infty,\infty}(\R^2)} \|g\|_{L^\infty(\R)} \geq 4 N$, and
such that~$g(0) = 0$, Theorem~\ref{th1.1} implies that there
is~$\e_N$ depending on~$u_0$ and on norms of~$f^h$ and~$g$ such
that~$u_0 + (f^h \otimes g,0)(x_1,x_2,\e_N x_3) $ belongs
to~${\mathcal G}$, where we have denoted~$f^h \otimes g(x) =
(f^1(x_h)g(x_3),f^2(x_h)g(x_3))$. Since~$\e_N$ only depends on norms
of~$f^h$ and~$g$, that implies that for any~$\lambda \in [-1,1]$,
the initial data~$u_0 + \lambda (f^h\otimes g,0)(x_1,x_2,\e_N x_3) $
also belongs to~${\mathcal G}$. Using~Proposition~\ref{sizedata} one
concludes that: passing through~$u_0,$ there exists  uncountable
number of segments of length~$N$ which are included in~${\mathcal
G}$. }
\end{rmk}
\begin{rmk}
\label{superposition} {\rm With the notation of Theorem~\ref{th1.1},
the data~$u_0(x)+(v^h_0 + \e w_0^h ,w_0^3)(x_1,x_2,\e x_3)$ belongs
to~${\mathcal G}$ as long as~$\e$ is small enough, so one can add to
that initial data any vector field of the type~$(v^{h (1)}_0 + \e_1
w_0^{h (1)} ,w_0^{3 (1)})(x_1,x_2,\e_1 x_3)$ and if~$\e_1$ is small
enough (depending on~$u_0$, on~$\e$, and on norms of~$v^h_0$, $w_0$,
$v^{h (1)}_0$ and~$w^{(1)}_0$), then the resulting vector field
belongs to~${\mathcal G}$. One thus immediately constructs by
induction superpositions of the type
$$
u_0(x)+ \sum_{j= 0}^J (v^{h (j)}_0 + \e_j w_0^{h (j)} ,w_0^{3
(j)})(x_1,x_2,\e_j x_3)
$$
which belong to~${\mathcal G}$ for small enough~$\e_j$'s, depending
on~$u_0$, on the norms of the profiles~$v^{h (j)}_0 $ and~$w^{(j)}_0
$, and on~$(\e_k)_{k < j}$.

Finally notice that one can also require the slow variation on the
profiles to hold on another coordinate than~$x_3$, up to obvious
modifications of the assumptions of the theorem.
}
\end{rmk}

\begin{rmk}
\label{rem  arkfrime}
{\rm In\ccite{cgp}, an even larger  initial data than the one of Theorem\refer{theoquasi2D}
 is constructed. However the size of the solution blows up when~$\e$ tends to $0$, and
 this is a strong obstacle to the use of a perturbative argument such as the one we will use here.
}\end{rmk}

\subsection{Scheme of the proof and organization of the paper}
Let us start by introducing some notation. We shall denote by~$C$ any constant, which may change from line to line, and we will write~$A \lesssim B$ if~$A \leq CB$.  In the following we shall
denote, for any point~$x = (x_1,x_2,x_3) \in\R^3$, its horizontal
coordinates by~$x_h \eqdefa (x_1,x_2)$. Similarly the horizontal
components of any vector field~$u=(u^1,u^2,u^3)$ will be denoted
by~$u^h \eqdefa(u^1,u^2)$ and the horizontal divergence will be
defined by~$\mbox{div}_h u^h \eqdefa \nabla^h \cdot u^h$,
where~$\nabla^h \eqdefa (\partial_1,\partial_2)$. Finally we shall
define the horizontal Laplacian by~$\Delta_h \eqdefa \partial_1^2 +
\partial_ 2^2$.
We shall often use the following shorthand notation for slowly
varying functions: for any function~$f$ defined on~$\R^3$, we write
\beq\label{notationslowvariation} [f]_\e (x_h,x_3) \eqdefa f(x_h, \e
x_3). \eeq In order to prove Theorem~\ref{th1.1}, we look for the
solution (which exists and is smooth for a short time depending
on~$\e$, due to classical existence theory) under the form \beq
\label{definreste} u_\e \eqdefa u_\e^{app}+R_\e \eeq where the
approximate solution~$u_\e^{app}$ is defined by the sum of the
global solution associated with~$u_0$ and the quasi-2D
approximation:
 \beq \label{definuapp}
 u^{app}_\e \eqdefa u+
[v^{(2D)}_\e]_\e \with v^{(2D)}_\e \eqdefa (v^h,0) + (\e
w^h_\e,w^3_\e)
\eeq
while
\begin{itemize}
\item
$u$ is the global smooth solution of~$(NS)$ associated with the
initial data~$u_0$;

\item
$v^h $ is the global smooth solution of the two dimensional
Navier-Stokes equation (with parameter~$y_3$ in~$\R$) with
pressure~$p_0$ and data~$v^h_{0}(\cdot, y_3)$
$$
(NS2D_3) \quad \quad \left\{
\begin{array}{c}
\p_t {v} ^h +v^h\cdot\na^h v^h-\Delta_h v^h = -\na^h p_0 \\
\dive_h\buh=0 \\
v^h |_{ t=0}=v^h_{0}(x_h,y_3) ;
\end{array}\right.
$$

\item $w_\e$ solves the linear equation with data~$w_0$ (and pressure~$
p_{\e,1}$)
$$
(T^\e_{v}) \left\{
\begin{array}{c}
\partial_tw_\e + v^h \cdot \nabla^h w_\e -\Delta_h w_\e
-\e^2\partial_3^2w_\e =
-( \nabla^h p_{\e,1} ,\e^2 \partial_{3}p_{\e,1} )\\
\dive w_\e = 0\\
{w_\e}|_{t=0} = w_0.
\end{array}
\right.
$$
\end{itemize}
We will also define the approximate pressure
 \beq \label{defpepsapp}
p_\e^{app} \eqdefa p+ [p_0+\e p_{\e,1}]_\e.
 \eeq
 The stability of
this approximate solution is described by the following proposition.
As in the rest of this
paper, we have used the following notation: if~$X$ (resp. $Y$) is a
function space over~$\R^2$ (resp. $\R$), then we write~$X_h$
for~$X(\R^2)$ and~$Y_v$ for~$Y(\R)$. We also denote the
space~$Y(\R;X(\R^2))$ by~$Y_v X_h$.
\begin{prop}\label{propreguluapp}
{\sl \label{stabilityuapp} For any positive~$\e_0$, the
family~$(u^{app}_\e)_{\e\leq \e_0}$ of approximate solutions is
uniformly bounded in~$ L^2(\R^+; L^\infty(\R^3)) $ and the
family~$(\nabla u^{app}_\e)_{\e\leq \e_0}$ is uniformly bounded in~$
L^2(\R^+;L^\infty_v(L^2_h))$.
}
\end{prop}
The size of the error term~$E_\e$ (this denomination will become
apparent in the next section) defined by
\beq
\label{definerrorterm}
E_\e \eqdefa (\partial_t-\D) u^{app}_\e + u^{app}_\e \cdot \nabla
u^{app}_\e + \nabla p^{app}_\e
\eeq
can be estimated as follows.
\begin{prop}
\label{estimaerroterm}
{\sl  The family~$(E_\e)_{\e\leq \e_0}$ of error
terms satisfies
$$
\lim_{\e \to 0} \|E_\e\|_{L^2(\R^+;\dot H^{-\frac 1 2})} =0.
$$
}
\end{prop}

\smallbreak The structure of this article is the following:

\begin{itemize}

\item the second section is devoted the proof of Theorem\refer{th1.1}
using the above two propositions;

\item the third section consists in proving
Proposition\refer{estimaerroterm} using estimates on the product in
anisotropic spaces;

\item we shall present the proof of some product laws in Sobolev
spaces in Appendix \ref{appendixa};

\item the proof of Proposition~\ref{propreguluapp} is  postponed to
Appendix~\ref{appendix}. Indeed most of the proof is actually contained
in Lemma 2.1 of~\cite{cg3}, apart from the fact that the global solution~$u$
satisfies the required properties. One way to avoid having to rely on that last
 result would be simply to replace, in the definition of~$ u^{app}_\e$,
 the solution~$u$ by any smooth approximation (that is possible due to the stability
 result of~\cite{GIP}). However we feel the result in itself is interesting
 so we prove in  Appendix~\ref{appendix} that any global solution associated to~$\dot H^\frac12(\R^3)$
 initial data belongs to~$ L^2(\R^+; L^\infty(\R^3)) $, and its gradient to~$
L^2(\R^+;L^\infty_v(L^2_h))$.

\end{itemize}

\setcounter{equation}{0}
\section{Proof of Theorem~\ref{th1.1}}\label{proofmainresult}
Assuming Proposition~\ref{estimaerroterm}, the proof of
Theorem\refer{th1.1} follows the sames lines as the proof of
Theorem~3 of\ccite{cg3}; we recall it for the reader's convenience.
Using the definition of the approximate
solution~$(u^{app}_\e,p^{app}_\e)$ given
in~(\ref{definuapp},\ref{defpepsapp}), and the error term~$E_\e$
given in~(\ref{definerrorterm}), we find that the remainder~$R_\e$
satisfies the following modified three-dimensional Navier-Stokes
equation
$$
(MNS_\e)\quad \quad \left\{
\begin{array}{c}
\p_tR_\e+R_\e\cdot\na R_\e-\Delta R_\e+u_\e^{app}\cdot\na
R_\e+R_\e\cdot\na u_\e^{app} =-E_\e -\na q_\e \\
\dive R_\e=0\andf R_\e|_{t=0}= 0,
\end{array}\right.
$$
with~$ q_\e \eqdefa p_\e - p_\e^{app}$. The proof of the theorem
reduces to the proof that~$(MNS_\e)$ is globally wellposed. We shall only write
the useful a priori estimates on~$ R_\e$, and leave to the reader the classical
arguments allowing to deduce the result. In particular we  omit the proof of the
 fact that the solution~$ R_\e$ constructed in this way is continuous in time with values in~$\dot H^\frac12(\R^3)$.

\smallskip

So let us define, for any~$\lambda >0$,
$$
R_\e^{\lam} (t) \eqdefa R_\e(t) \exp \left(-\lam \int_0^t V_\e(t')
\: dt'\right) \with V_\e(t) \eqdefa \| u_\e^{app}(t)\|_{L^\infty}^2
+ \| \na u_\e^{app}(t)\|_{L^\infty_vL^2_h}^2 .
$$
Note that Proposition~\ref{stabilityuapp} implies that~$
\displaystyle \int_{\R^+}V_\e(t) \: dt$ is uniformly bounded, by a
constant denoted by~$U$ in the following. Writing
also~$\displaystyle E_\e^{\lam} (t) \eqdefa E_\e(t) \exp \left(-\lam
\int_0^t V_\e(t') \: dt'\right) $, an~$\dot H^{\frac 1 2}$ energy
estimate on~$(MNS_\e)$ implies
$$
\longformule{ \frac12\frac d{dt} \|R_\e^{\lam} (t) \|_{\dot
H^\frac12}^2 + \|\na R_\e^{\lam} (t) \|_{\dot H^\frac12}^2 = - \lam
V_\e(t) \|R_\e^{\lam} (t) \|_{\dot H^\frac12}^2 -\bigl(E^\lam_\e
|R_\e^{\lam} \bigr)_{\dot H^\frac12} (t) } { {}- \Bigl( \exp \Bigl
(\lam \int_0^t V_\e(t')dt'\Bigr) R_\e^{\lam}\cdot \na R_\e^\lam+
u_\e^{app}\cdot \na R_\e ^{\lam} + R_\e ^{\lam}\cdot \na u_\e^{app}
\bigl | R_\e^{\lam} \Bigr)_{\dot H^\frac12} (t). }
$$
A law of product in Sobolev spaces (see \eqref{pdtlaw} in Appendix
\ref{appendixa}) and Proposition\refer{stabilityuapp} imply that
\beno \begin{split}
 \exp \Bigl (\lam \int_0^t V_\e(t')dt'\Bigr)
\bigl|(R_\e^{\lam}\cdot \na R_\e^\lam | R_\e^{\lam})_{\dot
H^\frac12}\bigr|\leq &Ce^{\lam U} \| R_\e^{\lam}(t)\|_{\dot H^1}^2
\|\na R_\e^{\lam}(t)\|_{\dot H^\frac12}\\
 \leq& Ce^{\lam U} \| R_\e^{\lam}(t)\|_{\dot
H^\frac12} \|\na R_\e^{\lam}(t)\|_{\dot H^\frac12}^2 .
\end{split}
\eeno
 Lemma 2.3 of\ccite{cg3} claims that \beq\label{lemma2.3}
\bigl|(b\cdot \na a + a\cdot \na b|b)_{\dot H^\frac12}\bigr| \leq C
\bigl( \|a\|_{L^\infty} + \|\na a\|_{L^\infty_v L^2_h}
\bigr)\|b\|_{\dot H^\frac12} \|\na b\|_{\dot H^\frac12},
 \eeq
so by
definition of~$V_\e$ we get
$$
\bigl|(u_\e^{app}\cdot \na R_\e^{\lam} + R_\e^{\lam}\cdot \na
u_\e^{app} | R_\e^{\lam} \bigr)_{\dot H^\frac12}\bigr| \leq \frac14
\|\na R_\e^{\lam}(t)\|_{\dot H^\frac12}^2 + CV_\e(t) \|
R_\e^{\lam}(t)\|_{\dot H^\frac12} ^2.
$$
Let us choose~$\lambda \geq C$. Then using the fact that
$$
\bigl|(E_\e^\lambda |R_\e^{\lam})_{\dot H^\frac12} (t)\bigr| \leq
\frac14 \|\na R_\e^{\lam}(t)\|_{\dot H^\frac12}^2 + C
\|E^\lambda_\e(t)\|_{\dot H^{-\frac12}}^2
$$
we obtain
$$
\frac d{dt} \|R_\e^{\lam} (t) \|_{\dot H^\frac12}^2 + \frac32 \|\na
R_\e^{\lam} (t) \|_{\dot H^\frac12}^2 \leq C \|E_\e(t)\|_{\dot
H^{-\frac12}}^2 + C e^{CU} \| R_\e^{\lam}(t)\|_{\dot H^\frac12}
\|\na R_\e^{\lam}(t)\|_{\dot H^\frac12}^2 .
$$
Since~${R_{\e}}_{|t=0}= 0$ and~$\displaystyle \lim_{\e \to 0}
\|E_\e\|_{L^2(\R^+;\dot H^{-\frac 1 2})} =0$ by
Proposition~\ref{estimaerroterm}, we deduce that as long
as~$\|R_\e^{\lam} (t) \|_{\dot H^\frac12} $ is smaller than~$ (4C
e^{ CU}) ^{-1}$, then for any~$\eta>0$ there is~$\e_0$ such that
$$
\forall \e \leq \e_0, \quad \|R_\e^{\lam} (t) \|_{\dot H^\frac12}^2
+ \frac34 \int_0^t \|\na R_\e^{\lam} (t') \|_{\dot H^\frac12}^2 \:
dt' \leq \eta,
$$
which in turn implies that
$$
\forall \e \leq \e_0, \quad\forall t \in \R^+, \quad \|R_\e^{\lam}
(t) \|_{\dot H^\frac12}^2 + \frac34 \int_0^t \|\na R_\e^{\lam} (t')
\|_{\dot H^\frac12}^2 \: dt' \leq \eta.
$$
That concludes the proof of the theorem. \qed

\section{The estimate of the error term}
In this section, we shall prove Proposition~\ref{estimaerroterm}.
Let us first remark that the error term~$E_\e$ can be decomposed as
$$
E_\e=E_\e^1+E_\e^2\with E_\e^2\eqdefa u\cdot \nabla
[v^{(2D)}_\e]_\e+ [v^{(2D)}_\e]_\e\cdot \nabla u.
$$
Thus the term~$E_\e^1$ is exactly the error term which appears
in\ccite{cg3}, and Lemmas 4.1, 4.2 and~4.3 of\ccite{cg3} imply that
\beq \label{estimrestcg3} \|E_\e^1\|_{L^2(\R^+;\dot H^{-\frac 1 2})}
\leq C_0\e^{\frac 1 3}. \eeq

In order to estimate the term~$E^2_\e$, let us first observe that
Lemmas 3.1 and 3.2 of\ccite{cg3} imply the following proposition.
\begin{prop}[\cite{cg3}]
\label{rappelcg3}
{\sl For any $s$ greater than~$-1$, for any~$\alpha \in
\N^3$ and for any positive~$t$, we have
$$
\|\partial^\al v^{(2D)}_\e(t )\|_{L^2_v\dot H^s_h} ^2+\int_0^t
\|\partial^\al \nabla^h v^{(2D)}_\e(t' )\|_{L^2_v\dot H^s_h} ^2 dt'
\leq C_0.
$$
}
\end{prop}
We shall also be using the following result, whose proof is
postponed to the end of this paragraph.
\begin{prop}
\label{smallatzero} {\sl The vector field~$v^{(2D)}_\e$ satisfies
\beq\label{estimateat0} \| v^{(2D)}_\e (
\cdot,0)\|_{L^\infty(\R^+;L^2_h)} + \| \nabla^h v^{(2D)}_\e (
\cdot,0)\|_{L^2(\R^+;L^2_h)} \leq C \e^\frac12. \eeq Furthermore,
$v^{(2D)}_\e$ is uniformly bounded in~$L^\infty(\R^+,\dot
H^\frac12(\R^3)) \cap L^2(\R^+,\dot H^\frac32(\R^3)) $. }
\end{prop}
Assuming this result, let us prove
Proposition~\ref{estimaerroterm}.\\

\no{\bf Proof of Proposition \ref{estimaerroterm}:} The stability
theorem of\ccite{GIP} claims in particular that
$$
\lim_{t\rightarrow \infty} \|u(t)\|_{\dot H^{\frac 1 2}} =0.
$$
As the set of smooth compactly supported divergence free vector
fields is dense in the space of~$\dot H^{\frac 1 2}(\R^3)$ divergence
free vector fields, this allows to construct for any
positive~$\eta$, a family~$ (t_j)_{1\leq j \leq N} $ of positive
real numbers and a family~$(\phi_j)_{1\leq j\leq N}$ of smooth
compactly supported divergence free vector fields such that
(with~$t_0=0$) \beq \label{approximcorGIPelem} \|\underline u_\eta
\|_{L^\infty(\R^+;\dot H^{\frac 1 2})} \leq \eta \with \underline
u_\eta \eqdefa u-\widetilde u_\eta\andf \widetilde u_\eta
(t,x)\eqdefa \sum_{j=1} ^N {\bf 1}_{[t_{j-1},t_j]}(t)\phi_j(x). \eeq
Then, for any positive~$\eta$, let us decompose~$E_\e^2$ as \beq
\label{decompresteGIP} E_\e^2= \underline E_{\e,\eta} +\wt
E_{\e,\eta} \with \underline E_{\e,\eta} \eqdefa \underline
u_\eta\cdot \nabla [v^{(2D)}_\e]_\e+ [v^{(2D)}_\e]_\e\cdot \nabla
\underline u_\eta. \eeq The term~$\underline E_{\e,\eta}$ will be
estimated thanks to the following lemma which is a generalization
of~(\ref{lemma2.3}).
\begin{lem}
\label{lemmaproduitalmostaniso}
{\sl Let~$a$ and~$b$ be two smooth
functions. We have
$$
\|ab\|_{\dot H^{\frac 1 2}} \leq C\|a\|_{\dot H^{\frac 12 }}\bigl(
\|\nabla^h b \|_{L^\infty_v (L^2_h)} +\|b\|_{L^\infty}
+\|\partial_3b\|_{L^2_v( \dot H^{\frac 1 2}_h)}\bigr).
$$
}
\end{lem}
\begin{proof}
For any function~$f$ in~$\dot H^{\frac 1
2} (\R^3)$, one has
\beq \label{lemmaproduitalmostanisodemoeq1} \|f\|_{\dot
H^{\frac 12 }} \leq \|f\|_{L^2_h\dot H_v^{\frac 12 }}+
\|f\|_{L^2_v\dot H_h^{\frac 12 }} .
 \eeq
 That  estimate may be proved simply by Plancherel's theorem
  (see for instance the end of the proof of Lemma~2.3 of~\cite{cg3}).

 Now we observe that by
two-dimensional product laws (taking $s=\frac12$ and $d=2$ in
\eqref{pdtlaws} of   Appendix \ref{appendixa}), one has for
any~$x_3$ in~$\R$
$$
\|a(\cdot,x_3)b(\cdot,x_3)\|_{\dot H^{\frac 1 2}_h} \leq C
\bigl(\|a(\cdot,x_3)\|_{_{\dot H^{\frac 1 2}_h}}\|\nabla^h
b(\cdot,x_3)\|_{L^2_h} +\|a(\cdot,x_3)\|_{\dot H^{\frac 1 2 }_h}
\|b(\cdot,x_3)\|_{L^\infty_h} \bigr).
$$
One has of course
\beq
\label{inegdebileaniso}
s\leq 0
\Longrightarrow \|a\|_{\dot H^s} \leq \|a\|_{L^2_v(\dot H^s_h)}
\andf s\geq 0\Longrightarrow \|a\|_{L^2_v(\dot H^s_h)} \leq
\|a\|_{\dot H^s}
\eeq
so taking~$s = 1/2$ gives \ben
\|ab\|_{L^2_v\dot H^{\frac 1 2}_h} & \leq & C \|a\|_{L^2_v(\dot
H^{\frac 1 2}_h)} \bigl( \|\nabla^h b \|_{L^\infty_v (L^2_h)}
+\|b\|_{L^\infty}\bigr)\nonumber\\
\label{lemmaproduitalmostanisodemoeq2} & \leq & \|a\|_{\dot H^{\frac
1 2}} \bigl( \|\nabla^h b \|_{L^\infty_v (L^2_h)}
+\|b\|_{L^\infty}\bigr). \een Now let us estimate~$\|ab\|_{L^2_h\dot
H^ {\frac 1 2}_v}$. A law of product in the vertical variable
(taking $s=\frac12$ and~$d=1$ in~\eqref{pdtlaws} of   Appendix
\ref{appendixa}) implies that for any~$x_h$ in~$\R^2$
$$
\|a(x_h,\cdot)b(x_h,\cdot)\|_{\dot H^{\frac 1 2}_v } \leq C \bigl (
\|a(x_h,\cdot)\|_{\dot H^{\frac 1 2}_v}
\|b(x_h,\cdot)\|_{L^\infty_v} + \|a(x_h,\cdot)\|_{L^2_v}\|\partial_3
b(x_h,\cdot)\|_{L^2_v}\bigr).
$$
Taking the $L^2$ norm in the horizontal variable gives
$$
\|ab\|_{L^2_h\dot H^{\frac 1 2}_v} \leq C \bigl( \|a\|_{L^2_h\dot
H^{\frac 1 2}_v} \|b\|_{L^\infty}
+\|a\|_{L^4_hL^2_v}\|\partial_3b\|_{L^4_hL^2_v} \bigr).
$$
Using Minkowski's inequality, we get that
$$
\|ab\|_{L^2_h\dot H^{\frac 1 2}_v} \leq C \bigl( \|a\|_{L^2_h\dot
H^{\frac 1 2}_v} \|b\|_{L^\infty}
+\|a\|_{L^2_vL^4_h}\|\partial_3b\|_{L^2_vL^4_h}\bigr).
$$
Then using the Sobolev embedding~$\dot H^{\frac
1 2}_h \hookrightarrow L^4_h$ and\refeq{inegdebileaniso}, we infer
 \beno
\|ab\|_{L^2_h\dot H^{\frac 1 2}_v} & \leq & C \bigl(
\|a\|_{L^2_h\dot H^{\frac 1 2}_v} \|b\|_{L^\infty} +
\|a\|_{L^2_v\dot H^{\frac 1 2}_h} \|\partial_3 b\|_{L^2_v\dot
H^{\frac
1 2}_h} \bigr)\\
& \leq & C \|a\|_{\dot H^{\frac 1 2}} \bigl( \|b\|_{L^\infty}
+\|\partial_3 b\|_{L^2_v\dot H^{\frac 1 2}_h}\bigr). \eeno Together
with\refeq{lemmaproduitalmostanisodemoeq1}
and\refeq{lemmaproduitalmostanisodemoeq2}, this proves
Lemma~\ref{lemmaproduitalmostaniso}.
\end{proof}
That lemma allows
to obtain the required estimate for~$\underline E_{\e,\eta} $. Using
the divergence free condition, we indeed have that
$$
\underline E_{\e,\eta} = \dive (\underline u_\eta\otimes
[v^{(2D)}_\e]_\e+ [v^{(2D)}_\e]_\e\otimes \underline u_\eta).
$$
So the above lemma implies that for any positive time~$t$
 \beno
\|\underline E_{\e,\eta} (t)\|_{\dot H^{-\frac 1 2}} & \leq & C
\|\underline u_\eta \otimes [v^{(2D)}_\e]_\e + [v^{(2D)}_\e]_\e
\otimes
\underline u_\eta \|_{\dot H^{\frac 1 2}} (t)\\
& \leq & C \|\underline u_\eta\|_{\dot H^{\frac 1 2}} \bigl(
\|\nabla^h [v^{(2D)}_\e]_\e\|_{L^\infty_vL^2_h}
+\|[v^{(2D)}_\e]_\e\|_{L^\infty}+\|\partial_3[v^{(2D)}_\e]_\e\|_{L^2_v\dot
H^{\frac 1 2}_h}\bigr)(t). \eeno By definition of~$[\: \cdot \:
]_\e$ and using~(\ref{approximcorGIPelem}), we get
$$
\|\underline E_{\e,\eta}(t)\|_{\dot H^{-\frac 1 2}} \leq C\eta
\bigl(\| \nabla^h v^{(2D)}_\e\|_{L^\infty_vL^2_h}
+\|v^{(2D)}_\e\|_{L^\infty}+\e^{\frac 1 2}\|\partial_3
v^{(2D)}_\e\|_{L^2_v\dot H^{\frac 1 2}_h}\bigr).
$$
Proposition\refer{rappelcg3}, along with Proposition~\ref{propreguluapp},
 gives finally
 \beq
\label{prooftheofinaleq1} \|\underline E_{\e,\eta}\|_{L^2(\R^+;\dot
H^{-\frac 1 2})} \leq C_0\eta.
\eeq In order to estimate the term~$
\widetilde E_{\e,\eta} $ let us observe that thanks to the
divergence free condition, we have
\beq \label{widetildeE00} \wt
E_{\e,\eta} = \wt u^h_\eta \cdot \na^h [v^{(2D)}_\e]_\e+ \e\wt
u_\eta^3 [\partial_3 v^{(2D)}_\e]_\e + [v^{(2D)}_\e]_\e \cdot \nabla
\wt u_\eta.
 \eeq
 Using a 3D law of product (namely~\eqref{pdtlaw} in Appendix
 \ref{appendixa})
gives
$$
\|\e\wt u_\eta^3 [\partial_3 v^{(2D)}_\e]_\e\|_{\dot H^{-\frac 1 2}}  \leq  C\e \|\wt u_\eta \|_{\dot H^{\frac 1 2}}
\|[\partial_3 v^{(2D)}_\e]_\e\|_{\dot H^{\frac 12}} .
$$
This gives
\beq
\label{widetildeE01}
\|\e\wt u_\eta^3 [\partial_3 v^{(2D)}_\e]_\e\|_{\dot H^{-\frac 1 2}}  \leq   \e^{\frac 1 2} \|\wt u_\eta \|_{\dot H^{\frac 1 2}} \|v_\e^{(2D)}\|_{ \dot H^{\frac 3 2}}.
\eeq
The two other terms of\refeq{widetildeE00} are estimated using the following lemma.
\begin{lem}
\label{lemmainteracslowvarying}
{\sl Let~$a$ and~$b$ be two smooth
functions. We have
$$
\|ab\|_{\dot H^{-\frac 1 2}} \leq C\|a\|_{L^2_v\dot H^{\frac 1 2}_h}
\|b(\cdot ,0)\|_{L^2_h} + C\|x_3 a \|_{L^2}
\|\partial_3b\|_{L^\infty_v\dot H^{\frac 1 2}_h}.
$$
}
\end{lem}
\begin{proof}
Let us decompose~$b$ in the following way:
 \beq
\label{demointeractslowvareq1} b(x_h, x_3) = b(x_h,0)+\int_0^{x_3}
\partial_3b(x_h,y_3)dy_3.
 \eeq
  Laws of product for Sobolev spaces
on~$\R^2$ (see~\eqref{pdtlaw} in Appendix
 \ref{appendixa}) together with Assertion\refeq{inegdebileaniso} gives \ben
\|a (b_{|x_3=0}) \|_{\dot H^{-\frac 1 2} } & \leq & \|a (b_{|x_3=0})
\|_{L^2_v\dot H^{-\frac 1 2}_h} \nonumber\\
& \leq & \biggl( \int_\R \|a (\cdot, x_3) b(\cdot ,0)\|^2_{\dot
H^{-\frac 1 2}_h} dx_3\biggr) ^{\frac 1 2}\nonumber\\
& \leq & C\|b(\cdot, 0) \|_{L^2} \Bigl(\int_\R \|a(\cdot ,
x_3)\|_{\dot H^{\frac 1 2}_h} ^2 dx_3\Bigr)^{\frac12}\nonumber\\
\label{demointeractslowvareq2} & \leq & C \|a\|_{L^2_v\dot H^{\frac
1 2}_h} \|b(\cdot ,0)\|_{L^2_h}.
\een
In order to
use\refeq{demointeractslowvareq1}, let us observe that for
any~$x_3$, two-dimensional product laws give
\beno \biggl \| a(\cdot
,x_3) \int_0^{x_3} \partial_3b(\cdot, y_3) dy_3\biggl\|_{\dot H
^{-\frac 1 2}_h} & \leq & C \|a(\cdot, x_3)\|_{L^2_h} \biggl|
\int_0^{x_3 } \| \partial_3b (\cdot
,y_3)\|_{\dot H^{\frac 1 2}_h } dy_3\biggr|\\
& \leq & C |x_3| \|a(\cdot ,x_3)\|_{L_h^2} \|\partial_3
b\|_{L^\infty_v\dot H^{\frac 1 2}_h}. \eeno The above estimate
integrated in~$x_3$ together with\refeq{demointeractslowvareq1}
and\refeq{demointeractslowvareq2} gives the result.
\end{proof}

Now let us apply this lemma to estimate $\wt{u}_\eta^h\cdot\na^h
[v^{(2D)}_\e]_\e$  and~$[v^{(2D)}_\e]_\e\cdot\na\wt{u}_\eta$. We get
$$
\longformule{ \|\wt{u}_\eta^h\cdot\na^h [v^{(2D)}_\e]_\e(t)\|_{\dot
H^{-\frac 1 2}} \leq C\|\wt{u}_\eta^h(t,\cdot )\|_{L^2_v(\dot
H^{\frac 1 2}_h)} \|\na^h v^{(2D)}_\e(t,\cdot ,0)\|_{L^2_h} } { {}+\e
\|x_3 \wt{u}_\eta^h(t)\|_{L^2} \|\partial_3\na^h
v_\e^{(2D)}(t,\cdot)\|_{L^\infty_v(\dot H_h^{\frac 1 2})} }
$$
and
$$
\longformule{ \|[v^{(2D)}_\e]_\e\cdot\na\wt{u}_\eta (t) \|_{\dot
H^{-\frac 1 2}} \leq C \|\na\wt{u}_\eta(t,\cdot)\|_{L^2_v(\dot
H^{\frac 1 2}_h)} \|v_\e^{(2D)}(t,\cdot,0)\|_{L^2_h} } { {}+\e
\|x_3\nabla \wt u_\eta(t,\cdot)\|_{L^2} \|\partial_3 v_\e^{(2D)}
(t,\cdot)\|_{L^\infty_v(\dot H^{\frac 1 2}_h)}. }
$$
By construction of~$\wt u_\eta$ and by Proposition
\ref{rappelcg3} and~\ref{smallatzero} (using the embedding of~$H^1(\R)$ into~$L^\infty(\R)$), together
with\refeq{widetildeE00} and\refeq{widetildeE01},  we
infer that \beq \label{widetildeE}
 \| \widetilde E_{\e,\eta}
\|_{L^2(\R^+;\dot H^{-\frac 1 2})} \leq C_\eta \e^\frac12 \eeq and
putting~(\ref{estimrestcg3}), (\ref{prooftheofinaleq1})
and~(\ref{widetildeE}) together proves
Proposition\refer{estimaerroterm}, up to the proof of
Proposition~\ref{smallatzero}.\qed

\medskip

Let us finally prove Proposition~\ref{smallatzero}.\\

\no{\bf Proof of Proposition~\ref{smallatzero}.}\  We recall
that~$v^{(2D)}_\e = (v^h,0) + (\e w^h_\e,w^3_\e)$, and due to the
form of~$(NS2D_3)$ it is clear that~$v^h(t,x_h,0) = 0$ for any~$(
t,x_h)$ in~$\R^+\times \R^2$. So it remains to estimate~$(\e w_\e^h,
w_\e^3) $. We first notice that due to Lemma 3.2 of~\cite{cg3},
$$
\| \e w_\e^h ( \cdot,0)\|_{L^\infty(\R^+;L^2_h)} + \| \e \nabla^h
w_\e^h ( \cdot,0)\|_{L^2(\R^+;L^2_h)} \leq C \e,
$$
so we are left with the computation of~$w_\e^3(t,\cdot,0)$. By
definition of~$w_\e $ we have
$$
\left\{
\begin{array}{c}
\partial_tw^3_\e + v^h \cdot \nabla^hw_\e^3 -\Delta_h w^3_\e = \e^2F_\e \\
{w^3_\e}_{|t=0} = w^3_0
\end{array}
\right. \with F_\e\eqdefa \partial_3^2w^3_\e - \partial_{3}p_{\e,1}.
$$
We shall start by writing an~$\dot H^{\frac 1 2}_h$ energy estimate
(with~$y_3$ seen as a parameter) which will imply
that~$w_\e^3(t,\cdot,0)$ is smaller than~$C\e$
in~$L^\infty(\R^+;\dot H^\frac12_h) \cap L^2(\R^+;\dot
H^\frac32_h)$. The result in the space~$L^\infty(\R^+;L^2_h) \cap
L^2(\R^+;\dot H^1_h)$ will follow by interpolation with a bound in a
negative order Sobolev space, given by Lemma 3.2 of\ccite{cg3}.

Let us start by the~$\dot H^{\frac 1 2}_h$ energy estimate. We claim
that there is a constant~$C_0$ such that for any~$\e \leq \e_0$,
\begin{equation}\label{estimateFeps}
\e\|F_\e \|_{L^2(\R^+;L^\infty_v\dot H^{-\frac 1 2}_h)} \leq C_0.
\end{equation}
Assuming~(\ref{estimateFeps}), an~$\dot H^{\frac 1 2}_h$ energy
estimate (joint with the fact that~$w^3_{\e |t=0} (\cdot ,0) = 0$)
gives directly that
$$
\|w^3_\e (\cdot ,0)\|_{L^\infty(\R^+;\dot H^\frac12_h)} + \|\nabla^h
w^3_\e (\cdot ,0)\|_{L^2(\R^+;\dot H^\frac12_h)} \leq C \e.
$$
But by  Lemma 3.2 of\ccite{cg3} we know that~$w^3_\e$ is uniformly
bounded, say in~$L^\infty(\R^+;L^\infty_v \dot H^{-\frac12}_h)$
and~$\nabla^h w^3_\e$ is uniformly bounded in~$L^2(\R^+;L^\infty_v
\dot H^{-\frac12}_h)$, so we get by interpolation that
$$
\| w_\e^3 ( \cdot,0)\|_{L^\infty(\R^+;L^2_h)} + \| \nabla^h w_\e^3 (
\cdot,0)\|_{L^2(\R^+;L^2_h)} \leq C \e^\frac12.
$$
This achieves \eqref{estimateat0}.

It remains to prove the claim~(\ref{estimateFeps}). On the one hand,
  Lemma 3.2 of~\cite{cg3} implies that
\beq\label{d32w3}
 \| \partial_3^2w^3_\e \|_{L^2(\R^+;L^\infty_v\dot
H^{-\frac 1 2}_h)}= \| \partial_3\na^h\cdot w^h_\e
\|_{L^2(\R^+;L^\infty_v\dot H^{-\frac 1 2}_h)} \leq C_0. 
\eeq The
estimate on the pressure seems slightly more delicate, but we notice
as in~\cite{cg3} that
 \beq\label{formulapressure} \displaystyle
-(\e^2 \partial_{3}^{2} + \Delta_{h}) p_{\e,1} = \dive_{h} \bigl(
v^h \cdot \nabla^{h}w^h_{\e} + \partial_{3} (w_\e^{3}v^h)\bigr).
\eeq 
Since~$ \e \partial_3 \dive_{h}(\e^2 \partial_{3}^{2} +
\Delta_{h})^{-1}$ is a uniformly bounded Fourier multiplier, this
implies by Sobolev embedding  that \beno
\begin{split}
\|\e \partial_3 p_{\e,1} \|_{L^2(\R^+;L^\infty_v\dot H^{-\frac 1
2}_h)} \lesssim& \|\e
\partial_3 p_{\e,1} \|_{L^2(\R^+;H^1_v\dot H^{-\frac 1 2}_h)}\\
\lesssim &\|v^h\cdot\nabla^{h}w^h_{\e}+\partial_{3}
w_\e^{3}v^h+w_\e^{3}\partial_{3} v^h\|_{L^2(\R^+;H^1_v\dot H^{-\frac
1 2}_h)}.
\end{split}
\eeno However thanks to \eqref{pdtlaw} and  using the estimates of
Lemmas 3.1 and 3.2 of \ccite{cg3}, one has \beno
\|v^h\cdot\nabla^{h}w^h_{\e}\|_{L^2(\R^+;H^1_v\dot H^{-\frac 1
2}_h)}\leq C\|v^h\|_{L^\infty(\R^+;H^1_v\dot H^{\frac 1
2}_h)}\|\nabla^{h}w^h_{\e}\|_{L^2(\R^+;H^1_vL^2_h)}\leq C. \eeno Due
to the divergence free condition of $w_\e,$ a similar estimate holds
for $\|\partial_{3} w_\e^{3}v^h\|_{L^2(\R^+;H^1_v\dot H^{-\frac 1
2}_h)}.$  While again thanks to \eqref{pdtlaw} and  using the
estimates of Lemmas 3.1 and 3.2 of \ccite{cg3}, we obtain \beno
\|w_\e^{3}\partial_{3} v^h\|_{L^2(\R^+;H^1_v\dot H^{-\frac 1
2}_h)}\leq C\|w_\e^{3}\|_{L^2(\R^+;H^1_v\dot H^{\frac 1
2}_h)}\|\partial_{3} v^h\|_{L^\infty(\R^+;H^1_vL^2_h)}\leq C. \eeno
As a consequence, we arrive at \beq \label{epsd3p} \|\e \partial_3
p_{\e,1} \|_{L^2(\R^+;L^\infty_v\dot H^{-\frac 1 2}_h)}\leq C_0.
\eeq The combination of~(\ref{d32w3}) and~(\ref{epsd3p}) proves the
claim, hence Estimate~(\ref{estimateat0}) of
Proposition~\ref{smallatzero}.

\medskip

Finally let us prove the bound in~$L^\infty(\R^+;\dot
H^\frac12(\R^3)) \cap L^2(\R^+,\dot H^\frac32(\R^3)) $
of~$v^{(2D)}_\e$. Actually the bound for~$(v^h, 0)$ follows from
Lemma 3.1 and Corollary 3.1 of~\cite{cg3}, so we just have to
concentrate on~$(\e w_\e^h,w_\e^3)$. Lemma~3.2  of~\cite{cg3} gives
that~$ w_\e $ is uniformly bounded in~$L^\infty(\R^+;\dot
H^\frac12(\R^3))$, as well as the fact that~$\nabla^h w_\e$ is
uniformly bounded in~$  L^2(\R^+;\dot H^\frac12(\R^3)) $ so by the
divergence free condition we only need to check that~$\e
\partial_3 w_\e^h$  is uniformly bounded in~$  L^2(\R^+;\dot H^\frac12(\R^3)) $. In fact, we
shall  prove first  that~$ \e \partial_3  w_\e^h$ is uniformly
bounded in~$ L^2(\R^+;L^2(\R^3)) $       and then that~$(\e
\partial_3)^2 w_\e^h$ is uniformly bounded in~$  L^2(\R^+;L^2(\R^3)) $, so
that the result will follow by interpolation, using Lemma~3.2
of~\cite{cg3} to deal with horizontal derivatives. Actually we shall
only concentrate on the first bound and leave the second to the
reader as it is very similar. Indeed we get by a standard ~$L^2$
energy estimate on~$w_\e^h$ that
$$
\frac12 \frac d {dt} \|w_\e^h\|_{L^2}^2 + \|\nabla^h w_\e^h\|_{L^2}^2 +  \|\e \partial_3w_\e^h\|_{L^2}^2 = - (v^h\cdot \nabla^h w_\e^h + \nabla^h p_{\e,1} | w_\e^h)_{L^2}.
$$
On the one hand we can write \beno | (v^h\cdot \nabla^h w_\e^h
|w_\e^h)_{L^2} |
  \leq   C \|v^h\|_{L^\infty}  \|\nabla^h w_\e^h\|_{L^2 }  \|   w_\e^h\|_{L^2 }
\eeno which implies that
$$
| (v^h\cdot \nabla^h w_\e^h   |w_\e^h)_{L^2} | \leq \frac14
\|\nabla^h w_\e^h\|_{L^2 }^2 + C   \|   w_\e^h\|_{L^2 }  ^2 \|v^h
\|_{L^\infty} ^2.
$$
To estimate the pressure term, we use again~(\ref{formulapressure})
which allows  to write (using the fact that~$\partial_3 w_\e^3 = -
\dive_h w_\e^h$)
\begin{eqnarray*}
\begin{split}
|( \nabla^h p_{\e,1} |w_\e^h)_{L^2} |&\leq C \int_\R
\|w_\e^h\|_{\dot H^\frac12_h} \Bigl( \| v^h\|_{\dot H^\frac12_h}  \|\nabla^h w_\e^h\|_{L^2_h}
+ \|   w_\e^3\|_{\dot H^\frac12_h}  \|   \partial_3v^h\|_{L^2_h} \Bigr) \: dx_3  \\
&\leq C  \|   w_\e^h\|_{L^2_v \dot H^\frac12_h} \Bigl( \|
v^h\|_{L^\infty_v\dot H^\frac12_h}  \|\nabla^h w_\e^h\|_{L^2 } + \|
w_\e^3\|_{L^2_v\dot H^\frac12_h}  \|
\partial_3v^h\|_{L^\infty_vL^2_h} \Bigr).
\end{split}
\end{eqnarray*}
This implies, after some interpolation estimates, that
$$
|( \nabla^h p_{\e,1} |w_\e^h)_{L^2} |\leq   \frac14
\|\nabla^h w_\e^h\|_{L^2 }^2 + C \|   w_\e^h\|_{L^2 }  ^2   \|
v^h\|_{L^\infty_v \dot H^\frac12_h} ^4   
  + C\|w_\e \|_{L^2_v\dot H^\frac12_h}^2 \|\partial_3 v
^h\|_{L^\infty_v L^2_h}.
$$
%
 Thus applying Gronwall's lemma ensures that
  \beno
\begin{split}
  \|w_\e^h(t)\|_{L^2}^2 &+ \int_0^t \|\nabla^h w_\e^h(t')\|_{L^2}^2 \: dt'
   + \int_0^t  \|\e \partial_3w_\e^h(t')\|_{L^2}^2  \: dt' \\
\leq& {} \bigl( {} C\|w_\e \|_{L^2(\R^+;L^2_v\dot H^\frac12_h)}^2
\|\partial_3 v^h\|_{L^\infty(\R^+;L^\infty_v L^2_h)} +
\|w_\e^h(0)\|_{L^2}^2
 \bigr) \\
 &\quad \quad  \times  \exp C \Bigl( \|  v^h\|_{L^4(\R^+; L^\infty_v \dot H^\frac12_h)} ^4+
   \|v^h\|_{L^2(\R^+; L^\infty)}^2
 \Bigr),
 \end{split}
  \eeno
so the results of Lemmas 3.1, 3.2 and Corollary 3.1 of~\cite{cg3}
allow to conclude that~$(\e \partial_3) w_\e^h$ is uniformly bounded
in~$ L^2(\R^+,L^2(\R^3)) $. The estimates are similar for~$(\e
\partial_3)^2 w_\e^h$, and that concludes the proof of the
proposition. \qed

\appendix
\setcounter{equation}{0}
\section{Product  laws in $\dot H^s(\R^d)$}\label{appendixa}

To prove the product laws in $\dot H^s(\R^d)$ as well as  Proposition
\ref{propregulu} below, we shall need some basic facts on
Littlewood-Paley analysis, which we shall recall here without proof
but refer for instance to~\cite{bcd} for all necessary details.
Let~$\widehat{\phi}$ (the Fourier transform of $\phi$) be a radial
function in~$\mathcal{D}(\mathbb{R}^{d})$ such
that~$\widehat\phi(\xi) = 1$ for~$|\xi|\leq 1$
and~$\widehat\phi(\xi) = 0$ for~$|\xi|>2$, and we define~$
\phi_{\ell}(x)= 2^{d \ell}\phi(2^{\ell}x).$ Then the frequency
localization operators are defined by
$$
S_{ \ell} = \phi_{\ell}\ast\cdot \quad\mbox{and}\quad \Delta_{\ell}
= S_{\ell +1} - S_{\ell}.
$$
Let~$f$ be in~$\mathcal{S}'(\mathbb{
  R}^{d})$, let~$p,q$ belong to~$[1,\infty]$, and let~$s < d/p$. We say that~$f$ belongs to~$\dot
B^{s}_{p,q}(\R^d)$ if and only if
\begin{itemize}
\item The partial sum $ \sum^{m}_{-m} \Delta_{\ell}f$ converges to $f$
  as a tempered distribution;
\item The sequence $\e_{\ell} = 2^{\ell s}\| \Delta_{\ell} f\|_{L^{p}}$
  belongs to $\ell ^{q}$.
\end{itemize}
We will also need a slight modification of those spaces, taking into
account the time variable; we refer to~\cite{CheLer} for the
introduction of that type of space  in the context of the
Navier--Stokes equations.
  Let $u(t,x)\in \mathcal{S}'(\mathbb{R}^{ 1 + d})$ and let $\Delta_ \ell $ be a
  frequency localization with respect to the~$x$ variable. We will say
  that $u\in \widetilde {L^{\rho}}(\R^+;\dot{B}^{s}_{p,q}(\R^d))$ if and only if
\begin{equation*}
   2^{\ell s}\|\Delta_ \ell u\|_{L^\rho(\R^+;L^{p} )} =\e_\ell \in \ell ^q,
\end{equation*}
and other requirements are the same as in the previous definition.
Note that there is an equivalent definition of Besov spaces in terms
of the heat flow: for any positive~$s$,
\[
\|u\|_{\dot B^{-s}_{p,r}} = \Bigl\|t^{\frac {s} 2} \|e^{t\Delta}u(t)
\|_{L^{p}}\Bigr\|_{L^{r}(\R^+,\frac {dt} t)}.
\]
%
 Now let us apply the above facts to study product laws in
$\dot H^s(\R^d)$. The proofs are very classical (see \cite{bcd} for
instance), and we present them here just for the readers'
convenience.

\begin{prop}\label{product}
{\sl  

$(i)   $ Let $a\in \dot{H}^{\frac{d-1}2}(\R^d)\cap
 \dot{H}^s(\R^d)$ and $b\in
 L^\infty(\R^d)\cap\dot{H}^{s+\frac12}(\R^d)$ for~$s>0.$  Then
 $ab\in \dot{H}^s(\R^d)$ and
 \beq\label{pdtlaws}
 \|ab\|_{\dot{H}^s}\lesssim
 \|a\|_{\dot{H}^{\frac{d-1}2}}\|b\|_{\dot{H}^{s+\frac12}}
 +\|a\|_{\dot{H}^s}\|b\|_{L^\infty}.
 \eeq

 $(ii)   $ Let $a\in \dot{H}^{s_1}(\R^d),$ $b\in \dot{H}^{s_2}(\R^d)$ with
$s_1+s_2>0$ and $s_1,s_2<\frac{d}2.$ Then $ab \in
\dot{B}^{s_1+s_2-\frac{d}2}_{2,1}(\R^d),$ and \beq \label{pdtlaw}
\|ab\|_{\dot{B}^{s_1+s_2-\frac{d}2}_{2,1}}\lesssim
\|a\|_{\dot{H}^{s_1}}\|b\|_{\dot{H}^{s_2}}. \eeq

 }
\end{prop}
\begin{proof} 
In what follows $(c_j)_{j\in\ZZ}$ (resp.
$(d_j)_{j\in\ZZ}$) will always be a generic element in the sphere of
$\ell^2$ (resp. $\ell^1$).

Thanks to Bony's decomposition\ccite{Bo81}, we
have~$\displaystyle ab =T_ab+T_ba+R(a,b) $, with
 $$  
T_ab=\sum_{j\in\ZZ}S_{j-1}a \D_jb \quad \mbox{and} \quad R(a,b)=\sum_{j\in\ZZ}\D_ja\wt{\D}_jb, \quad \mbox{while} \quad
\wt{\D}_jb=\sum_{\ell=-1}^1\D_{j+\ell}b.
$$

\no{\bf (i)}\ Bernstein's inequalities give
$$
\|S_ja\|_{L^\infty}\lesssim
c_j2^{\frac{j}2}\|a\|_{\dot{H}^{\frac{d-1}2}},
$$
so thanks to the support to the Fourier transform of
$T_ab$   we have 
\beno
\begin{split}
\bigl\|\D_\ell\bigl(T_ab\bigr)\bigr\|_{L^2}\lesssim
&\sum_{|j-\ell|\leq 5}\|S_{j-1}a\|_{L^\infty}\|\D_jb\|_{L^2}\\
\lesssim & \: c_\ell 2^{-\ell
s}\|a\|_{\dot{H}^{\frac{d-1}2}}\|b\|_{\dot{H}^{s+\frac12}}.
\end{split}
\eeno Similarly as $s>0,$ it follows that \beno
\begin{split}
\bigl\|\D_\ell\bigl(T_ba+R(a,b)\bigr)\bigr\|_{L^2}\lesssim
&\sum_{j\geq\ell-N_0}\|\D_ja\|_{L^2}\|S_{j+2}b\|_{L^\infty}\\
\lesssim &\sum_{j\geq\ell-N_0}c_j 2^{-j
s}\|a\|_{\dot{H}^{s}}\|b\|_{L^\infty}\lesssim c_\ell 2^{-\ell
s}\|a\|_{\dot{H}^{s}}\|b\|_{L^\infty}.
\end{split}
\eeno This achieves \eqref{pdtlaws}.

\no{\bf (ii)}\ The proof is similar to that of \eqref{pdtlaws}   noticing
that as $s_1<\frac{d}2,$ \beno
\begin{split}
\bigl\|\D_\ell\bigl(T_ab\bigr)\bigr\|_{L^2}\lesssim&\sum_{|j-\ell|\leq
5}\|S_{j-1}a\|_{L^\infty}\|\D_jb\|_{L^2}\\
\lesssim &\sum_{|j-\ell|\leq
5}c_j^22^{-j(s_1+s_2-\frac{d}2)}\|a\|_{\dot{H}^{s_1}}\|b\|_{\dot{H}^{s_2}}
\lesssim d_\ell
2^{-\ell(s_1+s_2-\frac{d}2)}\|a\|_{\dot{H}^{s_1}}\|b\|_{\dot{H}^{s_2}}.
\end{split}
\eeno 
The same estimate holds for
$\|\D_\ell\bigl(T_ba)\|_{L^2}.$
On the other hand, as $s_1+s_2>0,$ we deduce that \beno
\begin{split}
\bigl\|\D_\ell\bigl(R(a,b)\bigr)\bigr\|_{L^2}\lesssim& \sum_{j\geq
\ell-N_0}2^{\frac{d}2\ell}\|\D_ja\|_{L^2}\|\wt{\D}_jb\|_{L^2}\\
\lesssim& 2^{\frac{d}2\ell}\sum_{j\geq
\ell-N_0}c_j^22^{-j(s_1+s_2)}\|a\|_{\dot{H}^{s_1}}\|b\|_{\dot{H}^{s_2}}
\lesssim d_\ell
2^{-\ell(s_1+s_2-\frac{d}2)}\|a\|_{\dot{H}^{s_1}}\|b\|_{\dot{H}^{s_2}}.
\end{split}
\eeno This completes the proof of \eqref{pdtlaw}.
\end{proof}

\medskip

\section{Proof of Proposition~\ref{propreguluapp}}\label{appendix}

\setcounter{equation}{0}

Proposition~\ref{propreguluapp} follows from the next statement, as
the~$[v^{(2D)}_\e]_\e$ part was dealt with in Lemma 2.1
of~\cite{cg3}. It remians prove the next result.
\begin{prop}\label{propregulu}
{\sl \label{stabilityu} Let~$u_0 \in \dot H^\frac12(\R^3)$ be a
divergence free vector field generating a
  smooth, global solution~$u$ to $(NS).$ Then~$u$ belongs to~$ L^2(\R^+; L^\infty(\R^3)) $ and~$\nabla  u$ to~$
L^2(\R^+;L^\infty_v(L^2_h))$.
}
\end{prop}
\begin{proof}
 We shall start by proving that~$u$ belongs to the space~$ L^2(\R^+; L^\infty(\R^3)) $.
 Writing~$$u = e^{t\Delta} u_0 + w\quad\mbox{with}\quad w\eqdefa
- \int_0^t e^{(t-t') \Delta}{\mathbb  P} \mbox{div} \:(  u \otimes u) (t') \:
 dt',
$$ we only need to prove the result for~$w$ since by the continuous embedding of~$\dot H^\frac12(\R^3)$
into~$\dot B^{-1}_{\infty,2}(\R^3)$ it is immediate to check  using
the definition of Besov spaces via the heat flow, that~$e^{t\Delta}
u_0$ belongs to~$ L^2(\R^+; L^\infty(\R^3)) $. So let us concentrate
on~$w$. 
 By Theorems 1.1 and 2.1 of~\cite{GIP},
 $u$ belongs to~$\widetilde L^\infty (\R^+;\dot H^\frac12(\R^3))\cap \widetilde L^2(\R^+;\dot H^\frac32(\R^3)),$
 so we infer that $u$ belongs to $\widetilde {L^4}(\R^+;\dot H^1(\R^3))$ and
 therefore~$u \otimes u$ belongs to~$\widetilde L^2(\R^+,\dot B^\frac12_{2,1}(\R^3))$ thanks to \eqref{pdtlaw}.

 In particular there is a sequence~$d_\ell$
in the unit sphere of~$\ell^1_\ell(L^2_t)$ such that
\beq\label{estimatedeltaell}
 \|\Delta_\ell  (u \otimes u)(t)\|_{L^2} \lesssim d_\ell(t)
 2^{-\frac{\ell}2}\|u\|_{\widetilde L^\infty (\R^+;\dot H^\frac12)}\|u\|_{\widetilde L^2(\R^+;\dot H^\frac32)}.
\eeq 
By the Plancherel formula, we get
\begin{eqnarray*}
\|\Delta_\ell w (t) \|_{L^2} & \lesssim &  \int_0^t e^{-(t-t') 2^{2 \ell}} 2^\ell \|\Delta_\ell  (u \otimes u)(t')\|_{L^2} \: dt' \\ & \lesssim &  2^{\frac{\ell}2}\Bigl( \int_0^t e^{-(t-t') 2^{2 \ell}} d_\ell(t')\: dt'\Bigr)
\|u\|_{\widetilde L^\infty (\R^+;\dot H^\frac12)}\|u\|_{\widetilde L^2(\R^+;\dot H^\frac32)} . 
\end{eqnarray*}
Using the Cauchy-Schwarz inequality,  we infer
$$
\|\Delta_\ell w (t) \|_{L^2} \lesssim \|d_\ell (\cdot) \|_{L^2_t} 2^{-\frac{3\ell}2}
\|u\|_{\widetilde L^\infty (\R^+;\dot H^\frac12)}\|u\|_{\widetilde L^2(\R^+;\dot H^\frac32)} . 
$$
Then, using (anisotropic) Bernstein inequalities (see for instance\ccite{bcd}) we have
$$
\|\Delta_\ell w (t) \|_{L^\infty} + \|\nabla  \Delta_\ell w (t) \|_{L^2_h(L^\infty_v)} \lesssim 2 ^{\frac {3 \ell} 2} \|\Delta_\ell w (t) \|_{L^2}.
$$
Then we conclude that
$$
 \sum_\ell \bigl( \|\Delta_\ell w \|_{L^2(\R^+;L^\infty)} +  \|\Delta_\ell \nabla w \|_{L^2(\R^+;L^2_h(L^\infty))}\bigr)
 \lesssim \|u\|_{\widetilde L^\infty (\R^+;\dot H^\frac12)}\|u\|_{\widetilde L^2(\R^+;\dot
 H^\frac32)}.
 $$

Let us now prove the result for~$\nabla e^{t \Delta } u_0$. The  proof follows the lines of the equivalence of the dyadic and heat definitions of Besov spaces (see for instance\ccite{bcd}). Using Lemma~2.1 of\ccite{chemin20} and the Bernstein inequality, we get that
\begin{eqnarray*}
\|t^{\frac 1 2}  \nabla  e^{t\Delta}  u_0\|_{L^\infty_v(L^2_h)} & \leq & \sum_{j }
\|t^{\frac 1 2}2^{\frac {3j} 2}  \Delta_j e^{t\Delta}  u_0\|_{L^2}\\
 & \lesssim &  C \|u_0\|_{\dot H^{\frac 1 2}}   \sum_{j } t^{\frac 1 2} 2^{j}e^{-ct2^{2j}}c_{j}
\end{eqnarray*}
where~$(c_{j})_{j\in\ZZ}$ denotes, as in  all  this proof, a generic
element of the unit sphere of~$\ell^{2} $. Using that
\[
\sup_{t>0} \sum_{j} t^{\frac 1 2}   2^{j}e^{-ct2^{2j}} <\infty,
\]
 we infer, using the Cauchy-Schwarz  inequality (in~$j$) with the weight~$\ds 2^{j}e^{-ct2^{2j}}$,
 \begin{eqnarray*}
\|\nabla e^{t\D} u_0\|_{L^2(\R;L^\infty_v(L^2_h))} ^2 & = & \int_0^\infty
t \|\nabla e^{t\Delta}  u_0\|_{L^\infty_v(L^2_h)}^2\frac {dt} t \\
 &  \lesssim   & \|u_0\|^2_{\dot H^{\frac 1 2}}
\int_0^\infty \biggl( \sum_{j\in \ZZ}  t^{\frac 1 2}
2^{j}e^{-ct2^{2j}}c_{j}\biggr)^2\frac {dt} t \\
 &  \lesssim   & \|u_0\|^2_{\dot H^{\frac 1 2}}
 \!\int_0^\infty\! \biggl( \sum_{j\in \ZZ}  t^{\frac 1 2}
2^{j}e^{-ct2^{2j}}\biggr)\biggl( \sum_{j\in \ZZ}  t^{\frac 1 2}
2^{j}e^{-ct2^{2j}}c_{j}^2\biggr) \frac {dt} t\\
 &  \lesssim   & \|u_0\|^2_{\dot H^{\frac 1 2}}  \int_0^\infty \sum_{j\in \ZZ}  t^{\frac 1 2}
2^{j}e^{-ct2^{2j}}c_{j}^2 \frac {dt} t \,\cdotp
\end{eqnarray*}
Using Fubini's theorem, we infer
$$
\|\nabla e^{t\D} u_0\|_{L^2(\R;L^\infty_v(L^2_h))} ^2 \lesssim
\|u_0\|^2_{\dot H^{\frac 1 2}} \sum_{j\in \ZZ} c_{j}^2\int_0^\infty
t^\frac12 2^{j }e^{-ct2^{2j}}\frac {dt} t
$$
which gives the result.
\end{proof}

\bigskip

\noindent {\bf Acknowledgments.} Part of this work was done when P.
Zhang visited Laboratoire J.-L. Lions of Universit\'e Pierre et
Marie Curie, as well as when J.-Y. Chemin visited the Academy of
Mathematics $\&$ Systems Science, CAS. They would like to thank
these two institutions for their hospitality. I. Gallagher is
partially supported by the French Ministry of Research grant
ANR-08-BLAN-0301-01. P. Zhang is partially supported by NSF of China
under Grant 10421101 and 10931007, and the innovation grant from
Chinese Academy of Sciences under Grant GJHZ200829.

 \smallskip

 The authors are grateful to the anonymous referee
 for a careful reading of the manuscript and many useful suggestions.

\end{document}